\documentclass{article}
\usepackage{theorem}
\usepackage{amsmath,epsfig,graphicx,color}
\usepackage{mathrsfs,amssymb, amsfonts, latexsym, mathtools}

\newtheorem{lem}{Lemma}[section]

\newtheorem{thm}{Theorem}[section]
\newtheorem{cor}{Corollary}[section]
\newtheorem{remrem}{Remark}[section]

\begin{document}
\title{Asymptotics of densities of first passage times for spectrally negative L\'{e}vy processes}
\author{
Shunsuke Kaji\thanks{Department of Mathematics Meijo University, Tempaku, Nagoya 468-8502, Japan. kaji@gmath.meijo-u.ac.jp}
\quad and Muneya Matsui\thanks{Department of Business Administration, Nanzan University, 18
Yamazato-cho, Showa-ku, Nagoya 466-8673, Japan. e-mail:mmuneya@gmail.com}
}

\date{} 
\maketitle
We study a first passage time of a L\'{e}vy process over a positive constant level.
In the spectrally negative case we give conditions for absolutely continuity of the distributions of the first passage times.  
The tail asymptotics of their densities are also clarified, where the asymptotics 
depend on tail behaviour of the corresponding L\'evy measures. 
We apply our results to the mathematical finance, in particular, the credit default swap pricing. \vspace{2mm}\\
{\bf Key words}: First passage time; L\'{e}vy process; absolute continuity; spectrally negative \\
{\bf MSC2020 subject classifications}: 62E17, 60G51, 60G70.

\section{Introduction}

For a L\'{e}vy process $X=\{ X_t \}_{t\in[0,\infty)}$ starting from $X_0=0$ we consider a first passage time 
\begin{eqnarray*}
{\tau}_b=\inf \{ t>0| X_t>b \},\quad b>0 
\end{eqnarray*}
and a density of its distribution.

In insurance, the probability $P(\tau_b >t)$ and the density of $- \log P(\tau_b >t)$ are called the survival probability and the hazard function, respectively. In finance, through the paper Carr and Wu \cite{C-W} we refer to the target standardized credit contract as a unit recovery claim or URC. Before we do so, we assume that the asset price process of the firm follows a geometric L\'{e}vy process $e^{-X}$, which is a martingale, and notice that the last process drops below the debt $e^{-b}$ at default. This fundamental claim with expiry $T$ is simply a security paying off one dollar at default time $\tau_b$ if and only if default, whose event is $\{ \tau_b <T \}$, occurs. Thus, letting constant interest rate $r$, the value of URC is given by
\[ E[ e^{-r \tau_b} ; \tau_b <T]. \]
If the distribution $P(\tau_b \in \cdot)$ has a density and we know it, then we can clearly compute the last expectation.

On the other hand, the work on the distribution of the first passage time $\tau_b$ is a classical area of research. It is well-known that the tail distribution $P(\tau_b>t)$ is regularly varying  if and only if the L\'{e}vy process $X$ satisfies the Spitzer condition$($see Bertoin \cite{Ber}$)$. In addition, in the case that the $X$ is the spectrally negative L\'{e}vy process, we have obtained a remarkable identity involving the distributions of $X$ and $\tau_b$ $($see Bertoin \cite{Ber} and Sato \cite{S}$)$. In the other cases, the absolutely continuous properties of the distribution $P(\tau_b \in\cdot)$ are studied by Chaumont \cite{C}, Doney and Rivero \cite{D-R1}, and Savov and Winkel \cite{S-W}. However, we do not obtain a result on an analysis of a density function of a distribution of a first passage time for a L\'{e}vy process, except for the spectrally positive cases. In particular, Doney and Rivero \cite{D-R2} establish it for a subordinator and Peskir \cite{P} provides the explicit representation of the density of the first passage time for the stable process with no negative jumps.

In this note, for the spectrally negative L\'{e}vy process $X$ we compute the asymptotic behavior of the density function of the distribution $P(\tau_b \in dt)$ as $t$ goes to infinity. To do so, we need to use the central and the generalized central limit theorems and need to change the probability measure. In final, we apply our results for pricing the credit default swaps. By using Corollary $2.1$ we determine the asymptotic rate of the expectations
\[ E[ e^{-r \tau_b}] - E[ e^{-r \tau_b} ; \tau_b <T]~as~T\to \infty, \]
where the second term of the last expectations is  the value of URC.

\section{Notation and Main Results}

On a probability space $(\Omega,\mathcal{F},P)$ we define a spectrally negative L\'{e}vy process $X=\{ X_t \}_{t\in[0,\infty)}$ starting from $X_0=0$ and with the characteristic function
\[ E[e^{i\lambda X_1}]=e^{\Psi(i\lambda)},~\lambda\in\mathbf{R}, \]
where $\displaystyle \Psi(\xi)=\frac{\sigma^2}{2}\xi^2 +m \xi + \int_{(-\infty,0)} (e^{\xi z}-1-\xi z ) \nu(dz)$, $\xi\in\mathbf{C}$. Here, $\sigma\geq 0$, $m\in\mathbf{R}$, and the L\'{e}vy measure $\nu(dz)$ on $\mathbf{R}$ satisfies
\begin{eqnarray}
\label{condi:levy:measure}
\nu ( (0,\infty) ) = 0~and~\int_{(-\infty, 0)} \min \{ z^2, |z|\} \nu(dz) < \infty.    
\end{eqnarray}
Furthermore, we set a natural filtration $\{ {\mathcal{F}}_t^X \}_{t\in[0,\infty)}$ by the completions of the $\sigma$-algebras generated by $X$ and let ${\mathcal{F}}_\infty^X$ be the smallest $\sigma$-algebra included by $\bigcup_{t\ge 0}{\mathcal{F}}_t^X$. All martingales here are considered with respect to the filtered probability space $(\Omega,{\mathcal{F}}_\infty^X,\{ {\mathcal{F}}_t^X \}_{t\in[0,\infty)},P)$.

For the regular variation of a measure, we borrow the definition in Hult and Lindskog \cite{H-L}. Let $B_\epsilon:=\{ y \in \mathbf{R},\,|y|< \epsilon \}$ is the open ball centered at $0$ with radius $\epsilon$. Write $\mathbf{R}_0=\mathbf{R}\setminus \{0\}$ and let $\mathscr C_0$ denote the class of real valued bounded continuous function on $\mathbf{R}_0$ such that each $f\in \mathscr C_0$ 
vanishes on $B_\epsilon$ for some $\epsilon>0$. Let $M_0$ be the class of Borel measures on $\mathbf{R}_0$ whose restriction to $\mathbf{R} \setminus B_{r}$ is finite for each $r>0$.
 
We say that a measure $\nu \in M_0$ is regularly varying if there exists a nonzero $\mu\in M_0$ and a regularly varying function $r(t)$ such that 
\[
 r(t)\nu(t(\cdot)) \to \mu(\cdot)\quad \text{in}\ M_0\quad \text{as}\quad t \to \infty,
\]
where the convergence $\mu_n\to \mu$ in $M_0$ is in the sense of weak convergence, and is equivalent to $\int fd\mu_n \to \int f d\mu$ for all $f\in \mathscr C_0$ (see Theorem 3.1 $(\mathrm{ii})$ in Hult and Lindskog \cite{H-L}). 

In our case we take $\bar \nu(t)=\nu(t(-\infty,-1])$ with $t>0$, where $t(-\infty,-1]:= \{tx:x\in (-\infty,-1]\}$ (cf. Hult and Lindskog \cite{H-L} $(1.1)$), and assume that $\bar \nu(t)$ is regularly varying with index $-\alpha,\,\alpha>0$ (see definition in ${\bf B}_\alpha$, $1<\alpha \le 2$) so that for $r(t)=1/\bar\nu(t)$ and $x<0$
\[
 r(t)\nu(t(-\infty,x])= \frac{\nu((-x)t (-\infty,-1])}{\nu(t(-\infty,-1])} = \frac{\bar \nu((-x)t)}{\bar \nu(t)} \to (-x)^{-\alpha}
\]
as $t\to\infty$. This is the definition of regular variation for the measure $\nu$. Here, since we consider the spectrally negative L\'evy process, we take $x<0$. For later use, we also define 
$\bar \nu$ as a measure on $[0,\infty)$, i.e. for $x>0$
\[
 \bar \nu (x)= \bar \nu([x,\infty))=\nu((-\infty,-x]).
\]
Notice that $\bar \nu$ is of locally bounded variation on $(0,\infty)$. 
We denote by $BV_{loc}((0,\infty))$ a set of right-continuous functions of locally bounded variation on $(0,\infty)$.

We assume the condition $\mathbf{B}_\alpha$, $1<\alpha \le 2$: \\
$\nu((-\infty,-x])$ is left-continuous and 
\begin{eqnarray*}
\left\{
\begin{array}
[c]{l}%
\bar \nu(x) : regularly~varying~with~index~-\alpha~~~~~~~if~1<\alpha<2; \\
~~~~~~~~~~~~~~~~\displaystyle \int_{(-\infty, -1)} z^2 \nu(dz) < \infty~~~~~~~~~~~~~~~~~~if~\alpha=2.
\end{array}
\right.
\end{eqnarray*}
Notice that under the condition $\mathbf{B}_\alpha$, $1<\alpha \le 2$, $\bar \nu(x) \in BV_{loc}(0,\infty)$.

Before we obtain our main theorem involving a density of a first passage time
\begin{eqnarray*}
{\tau}_b=\left\{
\begin{array}
[c]{l}%
\inf \{ t>0| X_t>b \}~~~~~~~if~~~\{ \} \neq \emptyset;\\
~~~~~~~~~\infty~~~~~~~~~~~~~~~~~~if~~~\{ \} = \emptyset,
\end{array}
\right.
\end{eqnarray*}
over a positive constant level $b$, we need the following lemma, which directly follows from Theorem 1.5.12 in Bingham et al. \cite{B-G-T}. 

\begin{lem}
\label{lem:bgt:thm1512}
Assume the condition $\mathbf{B}_\alpha$, $1<\alpha<2$ and let $r(t) = 1 /\bar \nu(t)$. 
Then, there exists a regularly varying function $g(x)$ with index $1/\alpha$ such that $r(g(x))\sim x$ as $x\to\infty$. 
Here $g(x)$ is determined uniquely up to asymptotic equivalence\footnote{If $f(x) \sim g(x)$ as $x\to\infty$, then $f(x)$ and $g(x)$ are said to be asymptotic equivalence.}. 
\end{lem}

\begin{thm}
\label{thm:hit:stable:asmptotic}
Assume that $m=0$, $(1)$, $\mathbf{B}_\alpha$, $1<\alpha<2$, and for some $\beta:\alpha \le \beta <2$, 
\begin{eqnarray}
\label{thm:hit:stable:condi}
\liminf_{x \downarrow 0} \{ x^{\beta-2} \int_{(-x,0)} z^2 \nu(dz) \} >0. 
\end{eqnarray}
Then, the probability $P(\tau_b \in \cdot)$ is absolutely continuous on $(0,\infty)$ and its density function $p_b(t)$ satisfies
\begin{eqnarray*}
p_b(t) \sim b \cdot \frac{\sin (\frac{\pi}{\alpha})}{\pi} \Gamma(1+\frac{1}{\alpha}) (\frac{\alpha -1}{\Gamma(2-\alpha)})^{\frac{1}{\alpha}} \frac{1}{tg(t)}~as~t \to \infty,
\end{eqnarray*}
where the gamma function $\displaystyle \Gamma (x)=\int_0^\infty e^{-t} t^{x-1} dt$, $x>0$ and $g(t)$ is a function given by Lemma \ref{lem:bgt:thm1512}.
\end{thm}

\begin{thm}
\label{thm:hit:normal:asmptotic}
Assume that $m=0$, \eqref{condi:levy:measure}, $\mathbf{B}_2$, and
\begin{eqnarray}
\label{assump:thm:hit:normal:asmptotic}
\sigma >0 ~~ or ~~ \liminf_{x \downarrow 0} \{  x^{s-2} \int_{(-x,0)} z^2 \nu(dz) \} >0~for~some~s\in (0,2).
\end{eqnarray}
Then, the probability $P(\tau_b \in \cdot)$ is absolutely continuous on $(0,\infty)$ and its density function $p_b(t)$ satisfies
\begin{eqnarray*}
p_b(t) \sim \frac{b}{\sqrt{2\pi c^2}} t^{-\frac{3}{2}}~as~t \to \infty,
\end{eqnarray*}
where $\displaystyle c^2=\sigma^2 + \int_{(-\infty,0)} z^2 \nu(dz)$.
\end{thm}

To treat the case $m>0$, we consider the Esscher transform (see Section 3.3 of Kyprianou \cite{Ky} or Section 33 of Sato \cite{S}). 
For this we assume 
\begin{eqnarray}
\label{def:lambda:minus}
\lambda_- = \inf \{ \lambda<0 | \int_{(-\infty,-1)} e^{\lambda z} \nu(dz)<\infty \} < 0
\end{eqnarray}
and
\begin{eqnarray}
\label{assump:m:psi'}
m>0,\quad \Psi ' (\lambda_{-})<0,
\end{eqnarray}
where due to \eqref{def:lambda:minus} the cumulant function $\Psi (\lambda)$, $\lambda\in (\lambda_{-}, 0)$ is well-defined.
Since $\Psi''(\lambda)>0$ on $\lambda \in (\lambda_-,0)$ and $\Psi'(0)=m>0$, 
there exists a unique solution $\lambda_*\in (\lambda_-, 0)$ such that the equation $\Psi ' (\lambda_*)=0$, and hence, 
$\displaystyle \Psi (\lambda_*) = \min_{\lambda \in (\lambda_-, 0)} \Psi (\lambda)  <0$. 

\begin{cor}
\label{cor:asympt:pbt:general}
Assume that \eqref{condi:levy:measure}, \eqref{assump:thm:hit:normal:asmptotic}, \eqref{def:lambda:minus} and \eqref{assump:m:psi'}. 
Then, the probability $P(\tau_b \in \cdot)$ is absolutely continuous on $(0,\infty)$ and its density function $p_b(t)$ satisfies
\begin{eqnarray*}
p_b(t) \sim \frac{b}{\sqrt{2\pi d^2}} t^{-\frac{3}{2}} \exp \{ -\lambda_* b + \Psi (\lambda_*) t \}~as~t \to \infty,
\end{eqnarray*}
where $\displaystyle d^2=\sigma^2 + \int_{(-\infty, 0)} z^2 e^{\lambda_* z} \nu(dz)$.
\end{cor}

\begin{remrem}
In the case that $X$ has a continuous sample path as $\sigma=1$ ; $m\ge 0$, we have already known
\[ p_b(t) = \frac{b}{\sqrt{2\pi t^3}} e^{-\frac{(mt-b)^2}{2t}},~t>0 \]
$($see Examples $46.5$ and $46.6$ in Sato \cite{S}$)$. The tail behavior of this is generalized in new theorem and corollary.
\end{remrem}

\begin{remrem}
\label{rem:exp:minus:mu}
We assume \eqref{condi:levy:measure}, \eqref{assump:thm:hit:normal:asmptotic}, \eqref{def:lambda:minus} and \eqref{assump:m:psi'}. 
From Corollary \ref{cor:asympt:pbt:general} we can easily obtain that for any $\mu>0$
\[ E[e^{-\mu \tau_b} ; \tau_b \ge t] \sim \frac{b e^{-\lambda_* b}}{(\mu-\Psi (\lambda_*))\sqrt{2\pi d^2}}~t^{-\frac{3}{2}} e^{(-\mu+\Psi (\lambda_*)) t}~as~t\to\infty. \]
The last asymptotic estimate is applied for finance in the next section.
\end{remrem}

\section{Application to Finance}

Let $r$ and $T$ be positive constant interest rate and maturity, respectively. The asset price of the firm follows a geometric L\'{e}vy process $\{ S_t \}_{t \in [0,T]}$ satisfying
\[ S_t = e^{-X_t},~t \in [0,T] \]
with \eqref{condi:levy:measure} and \eqref{def:lambda:minus}. In mathematical finance, letting $\tilde{S}_t = e^{-rt} S_t$, $t \in [0,T]$, it is assumed that $P$ is the risk-neutral probability measure, that is, $\{ \tilde{S}_t \}_{t \in [0,T]}$ is a martingale under $P$. Here, we see that the last assumption is equivalent to
\[ m = -r + \frac{\sigma^2}{2} + \int_{(-\infty,0)} (e^{-z}-1+z) \nu(dz). \]

In the Carr and Wu \cite{C-W} they study a credit default swap, in particular, a unit recovery claim or URC, whose contract pays one dollar at default time $D$ if $D\le T$ and zero otherwise. In addition, in the Carr and Wu \cite{C-W} it is considered that the firm's value $\{ S_t \}_{t \in [0,T]}$ is bound below by a barrier $K$ before default, but drops below at default, where $K$ is the debt principal in $(0,1)$. Then, the well-known option pricing theory guarantees that $D=\tau_b$ as $b=-\log K$ and the value $U_T$ of URC is
\[ U_T = E[e^{-r\tau_b} ; \tau_b \le T]. \] 

Finally, to compute the last value, we add the conditions \eqref{assump:thm:hit:normal:asmptotic} and \eqref{assump:m:psi'}. 
By using the Laplace transform $\displaystyle E[e^{-r\tau_b}]=e^{-b \Psi^{-1}(r)}$ (cf. Theorem 46.3 in Sato \cite{S}) 
and the fact in Remark \ref{rem:exp:minus:mu}, we compute
\[ e^{-b \Psi^{-1}(r)} - U_T \sim \frac{b e^{-\lambda_* b}}{(r-\Psi (\lambda_*))\sqrt{2\pi d^2}}~T^{-\frac{3}{2}} e^{(-r+\Psi (\lambda_*))T}~as~T\to\infty, \]
where $\Psi^{-1}$ is the inverse function of $\Psi$.

\section{Proofs}

\subsection{Lemmas}

Before we begin to obtain our main results, we first need the following lemmas.

\begin{lem}
\label{lem:lim:levy:alpha}
If \eqref{condi:levy:measure}, $\mathbf{B}_\alpha$, $1<\alpha<2$, and \eqref{thm:hit:stable:condi} are valid, then
\[
 \lim_{x \to \infty} g^\leftarrow (x)x^{-2} \int_{(-x,0)} z^2 \nu (dz)=\frac{\alpha}{2-\alpha} , 
\]
where $g^{\leftarrow}(x)$ is the generalized inverse of $g$ (see (1.5.10) of Bingham et al. \cite{B-G-T}). 
\end{lem}
proof : 
We fist see the relation of $\bar \nu$ and $\nu$. For $x>y>0$ 
\[
 \nu((-\infty,-y))=\nu((-\infty,-y]) = \bar \nu ([y,\infty)) = \bar \nu ((y,\infty)) 
\] 
by left-continuity of $\nu$, so that 
\[
 \nu ((-x,-y)) = \bar \nu ((y,x)).
\]
Hence, we have
\[
 \int_{(-x,0)} z^2 \nu (dz) = \int_{(0,x)} -z^2 \nu(dz) 
\]
and thus by Karamata theorem for Stieltjes integral forms (for details, see Theorem 1.6.4 in Bingham et al. \cite{B-G-T}) 
\[
  \int_{(-x,0)} z^2 \nu (dz) \sim \frac{\alpha}{2-\alpha} x^2
\bar \nu (x ). 
\]
This implies 
\begin{align*}
  g^\leftarrow (x) x^{-2} \int_{(-x,0)} z^2 \nu (dz)
& \sim \frac{\alpha}{2-\alpha} g^\leftarrow (x) \bar \nu (x) \\
& \sim \frac{\alpha}{2-\alpha}  x \bar \nu(g(x)) \\
& \sim \frac{\alpha}{2-\alpha}, 
\end{align*}
where regular variation of $\bar \nu$ is used.

\begin{lem}
\label{lem:stable:bound:func}
Suppose \eqref{condi:levy:measure}, $\mathbf{B}_\alpha$, $1<\alpha<2$ and \eqref{thm:hit:stable:condi}.
Let $\delta\in (0,1/\alpha)$ and take $t_\alpha>0$ such that 
\begin{align}
 \label{condi:t:alpha}
\inf_{t\ge t_\alpha} \{ t/ g(t)^\delta\} >0\quad \text{and}\quad \inf_{t\ge t_\alpha} \{ t/g^\leftarrow (g(t))\} >0,
\end{align} 
which are possible by regular variation of both $g$ and $g^\leftarrow$ with indices $\alpha$ and $1/\alpha$, respectively. 
Then for this $t_\alpha$ there exists a positive constant $\kappa_{\alpha,\beta}$ such that
\[ |E[e^{i(\lambda/g(t)) X_t}]| \le e^{-\kappa_{\alpha,\beta}|\lambda|^\delta} \]
holds for any $|\lambda| >\pi $ and $t \ge t_\alpha$.
\end{lem}
proof: Recall $E[e^{i \lambda X_t}]=e^{t\Psi (i \lambda)}$ and observe
\[ |E[e^{i(\lambda/g(t)) X_t}]| \le e^{t \int_{(-\infty,0)} (\cos (\lambda z/g(t)) -1) \nu(dz)}. \]
To end this proof, it suffices to show that for some positive constant $\kappa_{\alpha,\beta}$
\begin{eqnarray}
\label{pf:thm:hit:stable:ineq}
t \int_{(-\infty,0)} (\cos (\lambda z/g(t)) -1) \nu(dz) \le -\kappa_{\alpha,\beta} |\lambda|^\delta 
\end{eqnarray} 
holds for all $|\lambda|>\pi$ and $t\ge t_\alpha$. 
To do so, we have
\begin{eqnarray*}
t~(\cos (\lambda z/g(t)) -1) &=& -2t \sin^2 ( \frac{\lambda z}{2g(t)} ) \\
                                        &\le& -\frac{2}{\pi^2} \lambda^2 \frac{t}{g(t)^2}~z^2~~~~~~~~~~~if~|\frac{\lambda z}{2 g(t)}| \le \frac{\pi}{2},~z<0,
\end{eqnarray*} 
where we have used the equality $\cos \theta -1=-2\sin^2(\frac{\theta}{2})$ and the inequality $\sin^2 \theta \ge \frac{4}{\pi^2} \theta^2$ on $|\theta| \le \frac{\pi}{2}$. Then, we have
\begin{eqnarray*}
t \int_{(-\infty,0)} (\cos (\lambda z/g(t)) -1) \nu(dz) &\le& 
-\frac{2}{\pi^2} \lambda^2 \frac{t}{g(t)^2} \int_{(-\frac{\pi g(t)}{|\lambda|},0)} z^2 \nu(dz) =:-2 I(\lambda). 
\end{eqnarray*}

We evaluate $I(\lambda)$, separating the cases by values of $\frac{\pi g(t)}{|\lambda|}$, i.e. 
consider intervals $(0,\epsilon_\beta),\,[\epsilon_\beta,T_\alpha]$ and $(T_\alpha,\infty)$ for 
$\frac{\pi g(t)}{|\lambda|}$ where constants $\epsilon_\beta<T_\alpha$ are defined in the following. 
For a positive constant $c_\beta$ take $\epsilon_\beta>0$ such that 
for any $x\in (0,\epsilon_\beta)$ 
\[
 x^{\beta-2}\int_{(-x,0)} z^2 \nu(dz) \ge c_\beta,
\] 
which is possible by the condition \eqref{thm:hit:stable:condi}. 
We require two conditions for $T_\alpha>0$. Take $T_\alpha$ such that for all $x>T_\alpha$
\[
 g^{\leftarrow}(x) x^{-2} \int_{(-x,0)} z^2 \nu (dz) \ge \frac{1}{2} \frac{\alpha}{2-\alpha}.
\]
This is assured by Lemma \ref{lem:lim:levy:alpha}. Moreover, fix $A>0$ and $\epsilon' >1/\alpha-\delta$ and then 
take $T_\alpha(A,\epsilon')$ so that it satisfies  
\[
 \frac{g^\leftarrow (y)}{g^\leftarrow (x)} \le A \max \big\{
(y/x)^{1/\alpha +\epsilon'},\,(y/x)^{1/\alpha -\epsilon'}\big\},\quad x,y\ge T_\alpha(A,\epsilon'). 
\]
This inequality is assured by Potter's theorem (\cite[p.25]{B-G-T}), since 
$g^{\leftarrow}$ is regularly varying with index $1/\alpha$. We put $T_\alpha =\max \{T_\alpha, \,T_\alpha(A,\epsilon') \}$.

Now we turn to the evaluation of $I(\lambda)$. \\
$(\mathrm{i})$ The case $0<\frac{\pi g(t)}{|\lambda|}<\epsilon_\beta$. By definition of $\epsilon_\beta$,
\begin{align*}
 I(\lambda) &= \big(
\frac{|\lambda|}{\pi g(t)}
\big)^\beta t \big( \frac{\pi g(t)}{|\lambda|}\big)^{\beta-2} \int_{(-\frac{\pi g(t)}{|\lambda|},0)} z^2 \nu(dz) \\
& \ge |\lambda|^\delta \big(
\frac{|\lambda|}{\pi g(t)}
\big)^{\beta-\delta} \frac{t}{(\pi g(t))^\delta} c_\beta \\
& \ge |\lambda|^\delta \frac{\epsilon^{\delta-\beta}_\beta}{\pi^\alpha} \frac{t}{g(t)^\delta} c_\beta. 
\end{align*}
Then, in view of we have for all $t\ge t_\alpha$ 
\[
 I(\lambda) \ge |\lambda|^{\delta} \kappa_{\alpha,\beta}' \quad \text{on}\ |\lambda| >\pi g(t)/\epsilon_\beta. 
\]

\noindent 
$(\mathrm{ii})$ The case $\epsilon_\beta \le \frac{\pi g(t)}{|\lambda|} \le T_\alpha$, we observe that 
\begin{align*}
 I(\lambda) &= t \big(\frac{|\lambda|}{\pi g(t)}\big)^\delta \big(
\frac{|\lambda|}{\pi g(t)}\big)^{2-\delta} \int_{(-\frac{\pi g(t)}{|\lambda|},0)} z^2 \nu(dz) \\
& \ge |\lambda|^\delta \frac{t}{(\pi g(t))^\delta} T_\alpha^{\delta-2} \int_{(-\epsilon_\beta,0)} z^2 \nu(dz) \\
& \ge |\lambda|^{\delta} \kappa_{\alpha,\beta}'', 
\end{align*}
where we use \eqref{condi:t:alpha} and the constant $\kappa_{\alpha,\beta}''$ is taken uniformly on $t\ge t_\alpha$. Thus 
\[
 I(\lambda) \ge |\lambda|^{\delta}\kappa_{\alpha,\beta}''\quad \text{on}\quad |\lambda| \le \pi g(t)/\epsilon_\beta\quad \text{or}
\quad |\lambda| \ge \pi g(t)/T_\alpha. 
\]

\noindent
$(\mathrm{iii})$ The case $T_\alpha <\frac{\pi g(t)}{|\lambda|}$. By definition of $T_\alpha$ we observe that 
\begin{align*}
 I(\lambda) &= \frac{t}{g^\leftarrow(\frac{\pi g(t)}{|\lambda|} )}g^\leftarrow(\frac{\pi g(t)}{|\lambda|} )
\big(\frac{\pi g(t)}{|\lambda|} \big)^{-2} \int_{(-\frac{\pi g(t)}{|\lambda|},0)} z^2 \nu(dz) \\
& \ge \frac{t}{g^\leftarrow(\frac{\pi g(t)}{|\lambda|} )} \frac{1}{2} \frac{\alpha}{2-\alpha}.
\end{align*}
Then again by definition of $T_\alpha$ and $|\lambda|>\pi$, 
\[
 \frac{g^\leftarrow(\frac{\pi g(t)}{|\lambda|} )}{g^\leftarrow (g(t))} \le A \big(
\frac{\pi}{|\lambda|}
\big)^{1/\alpha-\epsilon'}
\]
holds and thus 
\begin{align*}
 I(\lambda) &\ge |\lambda|^{-1/\alpha +\epsilon'} \frac{t}{A g^\leftarrow(g(t))} \frac{1}{2} \frac{\alpha}{2-\alpha} \\
& \ge |\lambda|^{-\delta} \kappa_\alpha,\quad \text{on}\quad |\lambda| < \pi g(t)/T_\alpha,
\end{align*}
where $\kappa_\alpha$ can be taken uniformly on $t\ge t_\alpha$.

Now put $\kappa_{\alpha,\beta}=2\min\{ \kappa_{\alpha,\beta}',\,\kappa_{\alpha,\beta}'',\,\kappa_\alpha\}$, we obtain 
the desired result.

\begin{lem}
\label{lem:limit:in:stable}
Under the conditions $m=0$, \eqref{condi:levy:measure}, and $\mathbf{B}_\alpha$, $1<\alpha<2$ we obtain
\[ \lim_{t\to\infty} E[e^{i \lambda \frac{X_t}{g(t)}}]=\varphi_\alpha(\lambda),~\lambda \in \mathbf{R}, \]
where $\displaystyle \varphi_\alpha(\lambda) \equiv \exp\Big\{ \int_{-\infty}^0 (e^{i \lambda z}-1-i\lambda z ) d (-z)^{-\alpha}  \Big\}.$
\end{lem}
proof: We use the decomposition: for $0<\epsilon<1$
\begin{eqnarray*}
 \log E[e^{i \lambda \frac{X_t}{g(t)}}] &=& -\frac{\sigma^2 \lambda^2}{2} \frac{t}{g(t)^2} +  
\Big( \int_{(-\infty,-\epsilon)} + \int_{[-\epsilon,0)} \Big) (e^{i\lambda z}-1-i\lambda z)~t~\nu(g(t)dz) \\
                                                  &=:& \mathrm{I}_1(t)+\mathrm{I}_2(t)+\mathrm{I}_3(t). 
\end{eqnarray*}
Since by Lemma \ref{lem:bgt:thm1512} $g(t)$ is regularly varying with index $1/\alpha>1/2$, we have
\begin{eqnarray}
\label{lim:I1}
\lim_{t\to \infty} t / g(t)^2= 0,
\end{eqnarray}
and thus $\displaystyle \lim_{t\to\infty} \mathrm{I}_1(t)=0$.

Before we compute $\mathrm{I}_2(t)$, recall that for $x<0$
\[
 t\nu(g(t)(-\infty,x]) \sim r(s)\nu(s(-\infty,x]) \to (-x)^{-\alpha}:=\mu((-\infty,x])\ \text{as}\ s=g(t)\to \infty,
\]
where Lemma \ref{lem:bgt:thm1512} is used in the first approximation. 
Write 
\begin{align*}
 \mathrm{I}_2(t) &= \int_{(-\infty,-\epsilon)}(e^{i\lambda z}-1-i\lambda z{\bf 1}_{\{-1\le z\}})~t~\nu(g(t)dz) \\
                    &\qquad -i\lambda \int_{(-\infty,-1)} z\,t \, \nu(g(t)dz)=:\mathrm{I}_{21}(t)+i\lambda \mathrm{I}_{22}(t).
\end{align*}
Notice that $z \mapsto (e^{i\lambda z}-1-i\lambda z{\bf 1}_{\{-1\le z\}})$ is bounded and continuous except for 
the discontinuity point $\{-1\}$ and $\mu(\{-1\})=\int_{\{-1\}} d(-x)^{-\alpha}=0$. 
Thus an extended version of weak convergence (e.g. Theorem 5.2 in Billingsley \cite{Billingsley:1968}) 
\[
 \mathrm{I}_{21}(t) \to \int_{(-\infty,-\epsilon)} (e^{i\lambda z}-1-i\lambda z{\bf 1}_{\{-1\le z\}}) d(-z)^{-\alpha}.
\]
Write 
\[
 \mathrm{I}_{22}(t) = \frac{t}{g(t)} \int_{(-\infty,-g(t))} (-z) \nu(dz) =\frac{-t}{g(t)} \int_{(g(t),\infty)} z\, \bar \nu(dz). 
\]
Recall that $\bar \nu \in BV_{loc}(0,\infty)$ and thus by the Karamata theorem for Stieltjes-integral forms(for details Theorem 1.6.5 in Bingham et al. \cite{B-G-T}), 
\[
 \mathrm{I}_{22}(t) \sim \frac{-\alpha}{1-\alpha}
t \nu(g(t)(-\infty,-1]) \to \frac{-\alpha}{1-\alpha} \int_{-\infty}^{-1} d (-z)^{-\alpha}=\frac{-\alpha}{1-\alpha},
\] 
which is equivalent to $\displaystyle \int_{-\infty}^{-1} (-z) d(-z)^{-\alpha}$, and thus
\begin{align}
\label{lim:I2}
\lim_{t\to\infty} \mathrm{I}_2(t) = \int_{(-\infty,-\epsilon)}(e^{i\lambda z}-1-i\lambda z) d(-z)^{-\alpha}. 
\end{align} 
Next we observe that 
\begin{align*}
 |\mathrm{I}_3(t)| &\le \int_{[-\epsilon,0)}|e^{i\lambda z}-1-i\lambda z| t \nu (g(t)dz) \\
&= \int_{[-g(t)\epsilon,0)} \Big|e^{i\lambda (z/g(t))}-1-i\lambda (z/g(t)) 
\Big| t \nu (dz) \\
& \le \frac{t \lambda^2}{g(t)^2} \Big(
\int_{[-1,0)} + \int_{[-g(t)\epsilon,-1)} 
\Big)z^2 \nu(dz) \\
&= \frac{t}{g(t)^2}\lambda^2 \int_{[-1,0)} z^2 \nu(dz) -\frac{t \lambda^2}{g(t)^2} \int_{(1,g(t)\epsilon]} z^2 \bar \nu(dz). 
\end{align*}
From \eqref{lim:I1} the first term of the right hand side of the observation converges to $0$ as $t\to\infty$. 
Again by the Karamata theorem $1.6.4$ in Bingham et al. \cite{B-G-T} 
\begin{align*}
 -\frac{\lambda^2 t}{g(t)^2} \int_{(1,g(t)\epsilon]}z^2 \bar \nu(dz) & \sim \frac{\lambda^2 t}{g(t)^2} \frac{\alpha}{2-\alpha}
(g(t)\epsilon)^2 \bar \nu (g(t)\epsilon) \\
&= \lambda^2\epsilon^2 \alpha/(2-\alpha) t \bar \nu(g(t)\epsilon) \\
&\sim \alpha/(2-\alpha) \lambda^2 \epsilon^{2-\alpha}\quad \text{as}\quad t\to\infty 
\end{align*}
and therefore, $\displaystyle \lim_{\epsilon\downarrow 0}\lim_{t\to\infty} \frac{\lambda^2 t}{g(t)^2} 
\int_{(1, g(t)\epsilon]} z^2 \bar \nu (d z)=0.$ Thus, $\displaystyle \lim_{\epsilon\downarrow 0}\lim_{t\to\infty} |\mathrm{I}_3(t)|=0.$ 
Now letting $\epsilon \downarrow 0$ in \eqref{lim:I2}, we obtain the desired result.

\begin{lem}
\label{lem:stable:constant}
Let $\varphi_\alpha(\lambda)$, $1<\alpha<2$ be a function given by Lemma \ref{lem:limit:in:stable}. Then,
\[ \int_{-\infty}^\infty \varphi_\alpha(\lambda) d\lambda =2~\Gamma(1+\frac{1}{\alpha}) (\frac{\alpha -1}{\Gamma(2-\alpha)})^{\frac{1}{\alpha}}  \sin (\frac{\pi}{\alpha}). \]
\end{lem}
proof: 
First, we check
\begin{eqnarray}
\label{pf:lem:stable:constant:1}
\varphi_\alpha(\lambda) = \exp\Big\{ -c_\alpha |\lambda|^\alpha ( 1 + sgn \lambda \cdot i \tan\frac{\pi\alpha}{2} ) \Big\},
\end{eqnarray} 
where $c_\alpha=\frac{\Gamma(2-\alpha)}{1- \alpha}\cos (\frac{\pi\alpha}{2})>0$. From Eq. (14.19) of Sato \cite{S}, we have 
\begin{align*}
 \int_{-\infty}^0 (e^{i\lambda z}-1-i\lambda z) d(-z)^{-\alpha} &= \alpha 
|\lambda|^\alpha \int_0^\infty (e^{i r(-sgn \lambda)}-1-i r(-sgn \lambda))r^{-\alpha-1} d r   \\
& =\alpha |\lambda|^{\alpha}\Gamma(-\alpha) e^{i (sgn \lambda) \pi \alpha/2} \\
& = - |\lambda|^\alpha c_\alpha (1+i \tan(\frac{\pi \alpha}{2})\cdot sgn \lambda),
\end{align*}
where in the first line we changed variables. Hence, we obtain \eqref{pf:lem:stable:constant:1}.

Next, by \eqref{pf:lem:stable:constant:1} we observe
\begin{eqnarray*}
\int_{-\infty}^\infty \varphi_\alpha(\lambda) d\lambda 
&=& \int_0^\infty e^{-c_\alpha \lambda^\alpha ( 1 + i \tan(\frac{\pi\alpha}{2}) )} d\lambda + 
\int_0^\infty e^{-c_\alpha \lambda^\alpha ( 1 - i \tan(\frac{\pi\alpha}{2}) )} d\lambda \\                    
&=& \alpha^{-1} {c_\alpha}^{-\frac{1}{\alpha}}  \int_0^\infty 
\big(e^{- i \tan(\frac{\pi\alpha}{2}) \cdot \eta}
+e^{i \tan(\frac{\pi\alpha}{2}) \cdot \eta}
\big)~e^{-\eta} \eta^{\frac{1}{\alpha}-1} d\eta \\
&=& \alpha^{-1} {c_\alpha}^{-\frac{1}{\alpha}} ~\Gamma(\alpha^{-1})~
\big\{ (1 - i \tan(\frac{\pi\alpha}{2}))^{-\frac{1}{\alpha}} + ( 1 + i \tan(\frac{\pi\alpha}{2}))^{-\frac{1}{\alpha}}
\big\}, 
\end{eqnarray*}
where the last line of the observation holds by using, e.g. Example $2.15$ in Sato \cite{S}.

On the other hand, we compute
\[ ( 1 \pm i \tan(\frac{\pi\alpha}{2}))^{-\frac{1}{\alpha}} =  |\cos (\frac{\pi\alpha}{2})|^{\frac{1}{\alpha}}~e^{i(\pm \frac{\pi}{\alpha} \mp \frac{\pi}{2})}. \]
Hence, the observation and the computation yield
\begin{eqnarray*}
\int_{-\infty}^\infty \varphi_\alpha(\lambda) d\lambda 
&=& \alpha^{-1} {c_\alpha}^{-\frac{1}{\alpha}} ~\Gamma(\frac{1}{\alpha})~|\cos (\frac{\pi\alpha}{2})|^{\frac{1}{\alpha}}~
\big(e^{i(-\frac{\pi}{\alpha}+\frac{\pi}{2})} + e^{i(\frac{\pi}{\alpha}-\frac{\pi}{2})} \big) \\
&=& \frac{1}{\alpha} (\frac{\alpha -1}{\Gamma(2-\alpha)})^{\frac{1}{\alpha}} ~\Gamma(\frac{1}{\alpha}) \cdot 2\sin (\frac{\pi}{\alpha}),
\end{eqnarray*}
which provides the result.

\begin{lem}
\label{lem:gaussian:bounded:func}
If \eqref{condi:levy:measure}, $\mathbf{B}_2$, and \eqref{assump:thm:hit:normal:asmptotic} are valid, then there exists a positive constant $\kappa$ such that
\[ |E[e^{i(\frac{\lambda}{\sqrt{t}}) X_t}]| \le e^{-\kappa |\lambda|^s} \]
holds for any $(\lambda,t) \in D$, where $s$ is given by \eqref{assump:thm:hit:normal:asmptotic}
 and $D=\{(\lambda,t)\in\mathbf{R}\times(0,\infty) | |\lambda|>1, t>1 \}$.
\end{lem}
proof: Since we see $E[e^{i \lambda X_t}]=e^{t\Psi (i \lambda)}$, we have
\[ |E[e^{i(\frac{\lambda}{\sqrt{t}}) X_t}]| = e^{-\frac{\sigma^2}{2}\lambda^2} \cdot e^{t \int_{(-\infty,0)} (\cos \frac{\lambda z}{\sqrt{t}} -1) \nu(dz)}. \]
In the case $\sigma>0$ the last equality yields the desired result. Thus, to end this proof, in the cases that $\sigma=0$ and there are $\kappa_0$, $\delta>0$ such that 
for all $0<r<\delta$
\begin{eqnarray}
\label{pf:lem:gaussian:bounded:func:ineq1}
\int_{(-r,0)} z^2 \nu(dz) / r^{2-s} > \kappa_0,
\end{eqnarray} 
which is guaranteed by \eqref{assump:thm:hit:normal:asmptotic}, it suffices to show that for some positive constant $\kappa$
\begin{eqnarray}
\label{pf:lem:gaussian:bounded:func:ineq2}
t \int_{(-\infty,0)} (\cos \frac{\lambda z}{\sqrt{t}} -1) \nu(dz) \le -\kappa |\lambda|^s,~(\lambda,t) \in D.
\end{eqnarray}

To do so, we have
\begin{eqnarray*}
t~(\cos \frac{\lambda z}{\sqrt{t}} -1) &=& -2t \sin^2 ( \frac{\lambda z}{2\sqrt{t}} ) \\
                                                    &\le& -\frac{2}{\pi^2} \lambda^2 z^2~~~~~~~~~~~~~~~if~|\frac{\lambda z}{2\sqrt{t}}| \le \frac{\pi}{2},~z<0,
\end{eqnarray*} 
where we have used the equality $\cos \theta -1=-2\sin^2(\frac{\theta}{2})$ and the inequality $\sin^2 \theta \ge \frac{4}{\pi^2} \theta^2$ on $|\theta| \le \frac{\pi}{2}$. Then,
\begin{eqnarray}
\label{pf:lem:gaussian:bounded:func:ineq3}
t \int_{(-\infty,0)} (\cos \frac{\lambda z}{\sqrt{t}} -1) \nu(dz) \le -\frac{2}{\pi^2} \lambda^2 \int_{(-\frac{\pi\sqrt{t}}{|\lambda|},0)} z^2 \nu(dz),~(\lambda,t) \in D.
\end{eqnarray}

On the other hand, since by \eqref{pf:lem:gaussian:bounded:func:ineq1} we have
\begin{eqnarray*}
\int_{(-\frac{\pi\sqrt{t}}{|\lambda|},0)} z^2 \nu(dz) \ge \kappa_0 (\frac{\pi\sqrt{t}}{|\lambda|})^{2-s}~~~~~~~if~0<\frac{\pi\sqrt{t}}{|\lambda|}<\delta,
\end{eqnarray*}
on $D_\delta \equiv \{(\lambda,t)\in D | 0<\frac{\pi\sqrt{t}}{|\lambda|}<\delta \}$ we have
\begin{eqnarray*}
\frac{2}{\pi^2} \lambda^2 \int_{(-\frac{\pi\sqrt{t}}{|\lambda|},0)} z^2 \nu(dz) \ge \frac{2}{\pi^s} \kappa_0 |\lambda|^s.
\end{eqnarray*}
In addition, on $D \setminus D_\delta$ we have
\begin{eqnarray*}
\frac{2}{\pi^2} \lambda^2 \int_{(-\frac{\pi\sqrt{t}}{|\lambda|},0)} z^2 \nu(dz) \ge \frac{2}{\pi^2} \int_{(-\delta,0)} z^2 \nu(dz)~|\lambda|^s,
\end{eqnarray*}
where due to $\mathbf{B}_2$ we notice $\displaystyle \int_{(-\delta,0)} z^2 \nu(dz)<\infty$. Thus, the last two inequalities and 
\eqref{pf:lem:gaussian:bounded:func:ineq3} provide \eqref{pf:lem:gaussian:bounded:func:ineq2}.
Hence, the proof is complete.

\begin{lem}
\label{lem:gaussian:limit}
Under the conditions $m=0$, \eqref{condi:levy:measure}, and $\mathbf{B}_2$ we have
\[ \lim_{t\to\infty} E[e^{i\lambda \cdot \frac{1}{\sqrt{t}} X_t} ] = e^{-\frac{c^2}{2} \lambda^2}, \lambda \in \mathbf{R}. \]
\end{lem}
proof : By using the inequality $|e^{ia}-1-ia-\frac{1}{2}(ia)^2| \le |a|^3 e^{|a|}$ we have
\[ \lim_{t\to\infty} t \{ e^{i(\frac{\lambda}{\sqrt{t}})z} - 1 - i(\frac{\lambda}{\sqrt{t}})z -\frac{1}{2} ( i(\frac{\lambda}{\sqrt{t}})z )^2 \} =0,~\lambda \in \mathbf{R},~z<0. \] 
In addition, by another inequality $|e^{ia}-1-ia-\frac{1}{2}(ia)^2| \le c_0 |a|^2$ for some positive constant $c_0$, we have that for all $t>0$, $\lambda \in \mathbf{R}$, and $z<0$
\[ | t \{ e^{i(\frac{\lambda}{\sqrt{t}})z} - 1 - i(\frac{\lambda}{\sqrt{t}})z -\frac{1}{2} ( i(\frac{\lambda}{\sqrt{t}})z )^2 \} | \le c_0 \lambda^2 z^2. \]
Thus, according to the dominated convergence theorem, the last convergence and inequality imply from $\mathbf{B}_2$ that
\begin{eqnarray}
\label{pf:lem:gaussian:limit:1}
\lim_{t\to\infty} t 
\int_{(-\infty,0)} \big\{ e^{i(\frac{\lambda}{\sqrt{t}})z} - 1 - i(\frac{\lambda}{\sqrt{t}})z -\frac{1}{2} ( i(\frac{\lambda}{\sqrt{t}})z )^2 
\big\} ~\nu(dz) = 0,~\lambda \in \mathbf{R}.
\end{eqnarray} 

On the other hand, if $t>0$, $\lambda\in \mathbf{R}$, then we have
\begin{eqnarray}
\label{pf:lem:gaussian:limit:2}
{} & & | E[e^{i\lambda \cdot \frac{1}{\sqrt{t}} X_t}] - e^{-\frac{c^2}{2} \lambda^2}|  \\
   &=& e^{-\frac{c^2}{2} \lambda^2} | \exp \Big \{ t \int_{(-\infty,0)} \big\{ 
e^{i(\frac{\lambda}{\sqrt{t}})z} - 1 - i(\frac{\lambda}{\sqrt{t}})z -\frac{1}{2} ( i(\frac{\lambda}{\sqrt{t}})z )^2\big \}~\nu(dz) \Big \} -1|, \nonumber
\end{eqnarray}
since we see $E[e^{i\lambda \cdot \frac{1}{\sqrt{t}} X_t}]=e^{t \Psi(i \cdot \frac{\lambda}{\sqrt{t}})}$.
Hence, \eqref{pf:lem:gaussian:limit:1} and \eqref{pf:lem:gaussian:limit:2} provide the desired result.

\subsection{Proofs of Theorems \ref{thm:hit:stable:asmptotic} and \ref{thm:hit:normal:asmptotic}}
We consider Theorem \ref{thm:hit:stable:asmptotic} (resp. \ref{thm:hit:normal:asmptotic}). According to Proposition $2.5$ (xii) in Sato \cite{S}, 
it follows from Lemma \ref{lem:stable:bound:func} (resp. \ref{lem:gaussian:bounded:func}) that for each $t>0$ a probability $P(X_t \in \cdot)$ is absolutely continuous on $\mathbf{R}$ and that its density function $p(t,b)$ is
\begin{eqnarray}
\label{def:ptb}
p(t,b) = \frac{1}{2\pi} \int_{-\infty}^\infty e^{-i b \lambda} E[e^{i \lambda X_t}] d\lambda,~t>0.
\end{eqnarray}
In addition, Theorem $46.4$ in Sato \cite{S} and Corollary $3$ in Bertoin \cite{Ber} $($see p$190)$ guarantee that there exists a density function $p_b(t)$ such that
\[ p_b(t) = b t^{-1} p(t,b),~t>0. \]
Therefore, to obtain the desired conclusions in Theorems \ref{thm:hit:stable:asmptotic} and \ref{thm:hit:normal:asmptotic}, it suffices to show that
\begin{eqnarray}
\label{ttimes:stable:asympt}
p(t,b) \sim \frac{\Gamma(1+\frac{1}{\alpha})}{\pi} (\frac{\alpha -1}{\Gamma(2-\alpha)})^{\frac{1}{\alpha}} \sin (\frac{\pi}{\alpha}) \frac{1}{g(t)}~as~t \to \infty
\end{eqnarray}
and
\begin{eqnarray}
\label{ttimes:gaussian:asympt}
p(t,b) \sim \frac{1}{\sqrt{2\pi c^2}} t^{-\frac{1}{2}}~~as~t \to \infty,
\end{eqnarray}
respectively.

To show \eqref{ttimes:stable:asympt}, by \eqref{def:ptb} we have
\begin{eqnarray*}
g(t)~p(t,b) = \frac{1}{2\pi} \int_{-\infty}^\infty e^{-i b (\lambda/g(t))} E[e^{i \lambda(X_t/g(t)) }] d\lambda,~t>0.
\end{eqnarray*}
Then, according to the dominated convergence theorem, Lemmas \ref{lem:stable:bound:func}, \ref{lem:limit:in:stable}, and \ref{lem:stable:constant} yield
\begin{eqnarray*}
\lim_{t\to\infty} g(t)~p(t,b) = \frac{1}{2\pi} \int_{-\infty}^\infty \varphi_\alpha(\lambda) d \lambda =\frac{\Gamma(1+\frac{1}{\alpha})}{\pi} (\frac{\alpha -1}{\Gamma(2-\alpha)})^{\frac{1}{\alpha}} \sin (\frac{\pi}{\alpha}).
\end{eqnarray*}
This means \eqref{ttimes:stable:asympt}.

Finally, we will show \eqref{ttimes:gaussian:asympt}. By \eqref{def:ptb} we have
\begin{eqnarray*}
\sqrt{2\pi t}~p(t,b) = \frac{1}{\sqrt{2\pi}} \int_{-\infty}^\infty e^{-i b (\frac{\lambda}{\sqrt{t}})} E[e^{i (\frac{\lambda}{\sqrt{t}}) X_t}] d\lambda,~t>0.
\end{eqnarray*}
Then, according to the dominated convergence theorem, Lemmas \ref{lem:gaussian:bounded:func} and \ref{lem:gaussian:limit} yield
\begin{eqnarray*}
\lim_{t\to\infty} \sqrt{2\pi t}~p(t,b) = \frac{1}{\sqrt{2\pi}} \int_{-\infty}^\infty e^{-\frac{c^2}{2} \lambda^2} d \lambda =\frac{1}{c},
\end{eqnarray*}
which means \eqref{ttimes:gaussian:asympt}.
Hence, the proof is complete. 

\subsection{Proof of Corollary $2.1$}
According to Theorem $33.2$ in Sato \cite{S} or Theorem 3.9 of Kyprianou \cite{Ky}, 
from \eqref{condi:levy:measure}, \eqref{def:lambda:minus}, and \eqref{assump:m:psi'} 
we obtain uniquely a probability measure $Q$ such that for all $t\ge 0$
\[ Q(A)=E[Z_t ; A],~A \in {\mathcal{F}}_t^X, \]
where $Z_t = e^{\lambda_* X_t -\Psi(\lambda_*) t}$, $t\in[0,\infty)$ is exponential martingale. 
Under the measure $Q$, $\{X_t\}_{t\in[0,\infty)}$ is again a L\'evy process such that 
\begin{eqnarray}
\label{measure:change:levy:process}
E_Q[e^{i\theta X_1}]=\exp\{ im^{\ast}\theta -\frac{\sigma^{\ast 2}}{2}\theta^2 
+ \int_{(-\infty,0)} (e^{i \theta z}-1-i \theta z ) \nu^\ast (dz) \},~\theta\in\mathbf{R},
\end{eqnarray} 
where $\sigma^{\ast}=\sigma$, $\nu^\ast(dz)= e^{\lambda_\ast x}\nu(dz)$ and 
\[
 m^\ast = m +\lambda_\ast \sigma^2 + \int_{(-\infty,0)}(e^{\lambda_\ast x}-1)x \nu(dx) = \Psi'(\lambda_\ast)=0,
\]
so that 
\[
 E_Q[e^{i\theta X_1}]=\exp\{-\frac{\sigma^2}{2}\theta^2 
+ \int_{(-\infty,0)} (e^{i \theta z}-1-i \theta z ) e^{\lambda_\ast z}\nu (dz) \}.
\]
Observe that 
\begin{eqnarray}
\label{measure:change:levy:measure:condi}
\int_{(-\infty,-1)} z^2 \nu^{\ast}(dz) = 
\int_{(-\infty,-1)} z^2 e^{\lambda_* z} \nu(dz) <\infty,
\end{eqnarray}
and then, for any $t \ge 0$
\[ P(A)=E_Q[1/Z_t ; A],~A \in {\mathcal{F}}_t^X \]
is valid and $\{ 1/Z_t \}_{t\in [0,\infty)}$ is the exponential martingale under $Q$. On the other hand, we clear see that 
\eqref{assump:thm:hit:normal:asmptotic} yields
\begin{eqnarray*}
\sigma >0 ~~ or ~~ \liminf_{x \downarrow 0} \{ \int_{(-x,0)} z^2 e^{\lambda_* z} \nu(dz) / x^{2-s} \} >0~for~some~s \in (0,2).
\end{eqnarray*}
From Theorem $2.2$ the last condition, \eqref{measure:change:levy:process}, and \eqref{measure:change:levy:measure:condi} imply that $Q(\tau_b \in \cdot)$ is absolutely continuous on $[0,\infty)$ and its density function $q_b(\cdot)$ satisfies
\begin{eqnarray}
\label{time:changed:qbt}
q_b(t) \sim \frac{b}{\sqrt{2\pi d^2}} t^{-\frac{3}{2}}~~as~t \to \infty.
\end{eqnarray}

On the set $\{ \tau_b \le \cdot \}\in {\mathcal{F}}_{\tau_b \wedge \cdot}^X$ \footnote{a $\wedge$ b=$\min$ $\{$a, b$\}$} we have $Z_{\tau_b \wedge \cdot} =Z_{\tau_b}$, and therefore, the optional sampling theorem implies that for all $t\ge 0$
\begin{eqnarray*}
P(\tau_b \le t) &=& E_Q[ 1/Z_t ; \tau_b \le t ] = E_Q[ E_Q[ 1/Z_t | {\mathcal{F}}_{\tau_b \wedge t}^X ] ; \tau_b \le t ] \\
                      &=& E_Q[ 1/Z_{\tau_b \wedge t} ; \tau_b \le t ]=E_Q[ 1/Z_{\tau_b} ; \tau_b \le t ] \\
                      &=& E_Q[ e^{-\lambda_* b +\Psi(\lambda_*)\tau_b} ; \tau_b \le t ] \\
                      &=& \int_0^t e^{-\lambda_* b +\Psi(\lambda_*)s} q_b(s) ds.
\end{eqnarray*} 
Thus, $P(\tau_b \in \cdot)$ is absolutely continuous on $[0,\infty)$ and its density function $p_b(\cdot)$ is
\begin{eqnarray*}
p_b(t) = e^{-\lambda_* b +\Psi(\lambda_*)t} q_b(t),~t\in[0,\infty).
\end{eqnarray*} 
Hence, the last equality and \eqref{time:changed:qbt} provide the desired conclusion.\\

\noindent {\bf Acknowledgments}
The earlier versions of contents in this paper have been presented at the annual
workshop ``Infinitely divisible processes and related topics'' held in Dec. 2022.
We acknowledge the comments and the hosts in the workshop. 
MM's research is partly supported by the JSPS Grant-in-Aid for Scientific Research C
(19K11868).

\end{document}